\documentclass[a4paper,11pt,twoside,reqno]{amsart}

\usepackage[utf8]{inputenc}
\usepackage[plainpages=false,pdfpagelabels=true]{hyperref}
\usepackage{amssymb,amsthm}
\usepackage[margin=1in]{geometry}
\usepackage[dvipsnames]{xcolor} 

\newtheorem{Theorem}{Theorem}[section]
\newtheorem{Prop}[Theorem]{Proposition}
\newtheorem{Lem}[Theorem]{Lemma}
\newtheorem{Thm}[Theorem]{Theorem}
\newtheorem{Cor}[Theorem]{Corollary}

\theoremstyle{definition}
\newtheorem{Dfn}[Theorem]{Definition}
\newtheorem{Bem}[Theorem]{Remark}

\newcommand{\tr}{\operatorname{Tr}}

\newcommand{\sff}{\mathrm{I\!I}}

\newcommand{\dv}{\text{ }dv_g}

\allowdisplaybreaks[1]

\newtheorem{innercustomgeneric}{\customgenericname}
\providecommand{\customgenericname}{}
\newcommand{\newcustomtheorem}[2]{%
  \newenvironment{#1}[1]
  {%
   \renewcommand\customgenericname{#2}%
   \renewcommand\theinnercustomgeneric{##1}%
   \innercustomgeneric
  }
  {\endinnercustomgeneric}
}

\newcustomtheorem{customthm}{Theorem}

\renewcommand{\epsilon}{\varepsilon}

\newcommand{\s}{\mathbb{S}}

\numberwithin{equation}{section}

\title{On the normal stability of triharmonic hypersurfaces in space forms}

\author{Volker Branding}
\address{University of Vienna, Faculty of Mathematics\\
Oskar-Morgenstern-Platz 1, 1090 Vienna, Austria}
\email{volker.branding@univie.ac.at}

\date{\today}



\thanks{The author gratefully acknowledges the support of the Austrian Science Fund (FWF) through the project ``Geometric Analysis of Biwave Maps'' (P34853).
}

\begin{document}
\begin{abstract}
This article is concerned with the stability of triharmonic maps and in particular
triharmonic hypersurfaces.
After deriving a number of general statements on the stability of triharmonic maps
we focus on the stability of triharmonic hypersurfaces in space forms,
where we pay special attention to their normal stability.

We show that triharmonic hypersurfaces of constant mean curvature in Euclidean space
are weakly stable with respect to normal variations
while triharmonic hypersurfaces of constant mean curvature in hyperbolic space are 
stable with respect to normal variations.

For the case of a spherical target we show that the normal index of the small proper triharmonic hypersphere
\(\phi\colon\s^m(1/\sqrt{3})\hookrightarrow\s^{m+1}\) is equal to one
and make some comments on the normal stability of the proper triharmonic Clifford torus.
\end{abstract}
\maketitle

\section{Introduction and results}
\textit{Harmonic maps} are a geometric variational problem with rich structure
that has many applications in geometry, analysis and theoretical physics.
The geometric setup is the following. We consider a smooth map \(\phi\colon M\to N\) between two Riemannian manifolds \((M,g)\) and \((N,h)\).
Then, the \textit{energy} of \(\phi\) is defined by
\begin{align}
\label{1-energy}
E(\phi)(=E_1(\phi))= \frac{1}{2}\int_M|d\phi|^2\dv\, .
\end{align}
Its critical points are governed by the vanishing of the so-called
\emph{tension field} $\tau(\phi)$, that is
\begin{align}
\label{harmonic-map-equation}
0=\tau(\phi):=\tr_g \bar\nabla d\phi,
\end{align}
where \(\bar\nabla\) represents the connection on the pull-back bundle \(\phi^\ast TN\).
Harmonic maps are precisely the solutions of equation \eqref{harmonic-map-equation}.

The harmonic map equation \eqref{harmonic-map-equation}
is a second order semilinear elliptic partial differential equation.
Due to its second order nature powerful tools such as the maximum principle help to obtain a deep understanding of both its analytic and geometric properties,
see for example \cite{MR2389639} for an overview on harmonic maps.

Another geometric variational problem which received growing attention in recent years is that of the so-called \emph{biharmonic maps}.
These maps are characterized as the critical points of the \emph{bienergy} for maps between two Riemannian manifolds, which is given by
\begin{align}
\label{2-energy}
E_2(\phi)=\frac{1}{2} \int_M|\tau(\phi)|^2\dv.
\end{align}
Here, the Euler-Lagrange equation is a fourth order semilinear elliptic 
partial differential equation
and is expressed by means of the vanishing
of the \emph{bitension field} $\tau_2(\phi)$, that is
\begin{align}
\label{biharmonic-map-equation}
0=\tau_2(\phi):=\bar\Delta\tau(\phi)+\tr_g R^N(d\phi(\cdot),\tau(\phi))d\phi(\cdot).
\end{align}
Here, \(\bar\Delta\) is the so-called \textit{rough Laplacian}, i.e., the connection Laplacian on \(\phi^\ast TN\). For more background and the current status of research on biharmonic maps we refer to \cite{MR2301373} and the recent book \cite{MR4265170}.

The fact that \eqref{biharmonic-map-equation} is of fourth order leads to a number of 
technical difficulties such that the analysis of biharmonic maps is not as complete as the one of harmonic maps.
A direct inspection of the Euler-Lagrange equation for biharmonic maps \eqref{biharmonic-map-equation}
shows that every harmonic map automatically solves the equation for biharmonic maps.
Due to this reason one of the major objectives in the analysis of biharmonic maps is to find non-harmonic 
solutions which are called \emph{proper biharmonic}. On the other hand, there are geometric configurations
in which every biharmonic map necessarily needs to be harmonic: If \(M\) is compact and \(N\) has non-positive
sectional curvature, then every biharmonic map needs to be harmonic, see \cite{MR4040175}
for an overview on further results of this kind.
For this reason one usually
considers a spherical target when constructing proper biharmonic maps.
Let us mention several results concerning the stability of biharmonic maps
which are closely connected to the main results of this paper.
The indices of various proper biharmonic maps to spheres were calculated in
\cite{MR2135286,MR2327029,MR4216418}
while the stability of biharmonic hypersurfaces was studied in \cite{MR4386842}.

A further higher order generalization of harmonic maps is the notion of \emph{triharmonic maps}.
Here, the starting point is the \emph{trienergy} of a map given by
\begin{align}
\label{3-energy}
E_3(\phi)=\frac{1}{2} \int_M|\bar\nabla\tau(\phi)|^2\dv.
\end{align}
The critical points of \eqref{3-energy} are characterized by the vanishing 
of the \emph{tritension field}
\begin{align}
\label{triharmonic}
0=\tau_3(\phi):=-\bar\Delta^2\tau(\phi)+ \tr_g R^N(\bar\nabla_{(\cdot)}\tau(\phi),\tau(\phi))d\phi(\cdot)
+\tr_g R^N(\bar\Delta\tau(\phi),d\phi(\cdot))d\phi(\cdot)
\end{align}
and are precisely triharmonic maps. 
The triharmonic map equation is a semilinear elliptic partial differential equation of order six,
the large number of derivatives gives rise to further technical difficulties in the mathematical analysis.
Again, a direct inspection of the Euler-Lagrange equation shows
that harmonic maps are always a solution of the equation for triharmonic maps \eqref{triharmonic}.
Hence, we are again interested in constructing solutions of \eqref{triharmonic}
which are non-harmonic and as in the case of biharmonic maps the latter are called
\emph{proper triharmonic}.
We would like to point out that a harmonic map
is always automatically also triharmonic, but the same is not true for a given biharmonic map, i.e.
a biharmonic map is not automatically triharmonic.

Let us give some motivation why, besides overcoming the technical difficulties, it should be of great
interest to investigate the properties of triharmonic maps.
While biharmonic maps seem to be
geometrically restricted in the sense that biharmonic curves need to have constant curvature
or that compact biharmonic hypersurfaces always seem to have constant mean curvature
(in the non-compact case there exist biharmonic hypersurfaces with non-constant mean curvature, 
see \cite{MR2975260})
these restrictions seem no longer in place in the triharmonic case:
Recently, in \cite{MR4308322} an explicit triharmonic curve
with non-constant geodesic curvature was constructed.

In the following, we will give a non-exhaustive overview on the current status of research on triharmonic maps.
Triharmonic curves of constant curvature in various ambient geometries
have been studied in \cite{MR4308322,MR4440045,MR4424860}, polyharmonic, and hence in particular
triharmonic, helices in space forms were classified in \cite{MR4542687}.
The first attempt of explicitly constructing triharmonic maps in the sphere was carried out in \cite{MR3403738} which was later systematically extended in \cite{MR3711937}.

Much effort has been paid to the case of triharmonic hypersurfaces.
Polyharmonic hypersurfaces in Riemannian space forms 
with constant norm of the shape operator were studied in \cite{MR4462636},
this analysis was recently extended to the pseudo-Riemannian case as well \cite{MR4552081}.
Triharmonic hypersurfaces in space forms which have constant mean curvature
and at most three distinct principal curvatures were classified in 
\cite{MR4444191}.
Very recently, much progress has been made regarding the classification of triharmonic hypersurfaces
in space forms \cite{MR4598081}. In particular, it was shown that any triharmonic hypersurface
of constant mean curvature in the hyperbolic space \(\mathbb{H}^{m+1}\) is actually minimal.
Triharmonic surfaces in homogeneous three manifolds (BCV spaces) have been studied in \cite{mor2023}.

Let us also give some examples of non-existence results for triharmonic maps.
Triharmonic isometric immersions from complete non-compact Riemannian manifolds 
into hyperbolic space have been investigated in \cite{MR3371364}
and it was shown that such immersions must be minimal whenever a certain energy is finite.
For triharmonic maps from Euclidean space to Riemannian manifolds a classification result
for finite energy solutions was given in
\cite[Section 4]{MR4007262}.
In \cite{MR4184658} a classification result for polyharmonic maps, which can also be applied to
the case of triharmonic maps, was established.

The current status of research on higher order variational
problems can be found in \cite{MR4106647}.

So far, the stability of triharmonic maps has not been investigated systematically.
The first article that took up this direction of research 
was published in Chinese \cite{MR1131891}.
According to the \emph{Mathematical Reviews} the author showed that a harmonic map
from a compact Riemannian manifold is a stable triharmonic map and also
that if one considers a proper biharmonic map which satisfies a conservation law
from a compact manifold into a manifold of positive sectional curvature,
then it is an unstable triharmonic map.
Besides that, the second variation formula for polyharmonic maps between Riemannian 
manifolds was derived in \cite{MR3007953} which of course includes the triharmonic case as well.

Now, let \(\phi\colon M\to N\) be a smooth triharmonic map.
The stability of a given triharmonic map is characterized by the second variation of the trienergy of a map \eqref{3-energy} evaluated at a critical point
which we denote by 
\(\operatorname{Hess}E_3(\phi)(V,W)\), where
\(V,W\in\Gamma(\phi^\ast TN)\).
We say that a triharmonic map \(\phi\) is \emph{stable} if
\begin{align*}
\operatorname{Hess}E_3(\phi)(V,V)>0\qquad \textrm{ for all } V\in\Gamma(\phi^\ast TN),\qquad V\neq 0
\end{align*}
and \emph{weakly stable} if
\begin{align*}
\operatorname{Hess}E_3(\phi)(V,V)\geq 0\qquad \textrm{ for all } V\in\Gamma(\phi^\ast TN),\qquad V\neq 0.
\end{align*}

Concerning the stability of arbitrary triharmonic maps we prove 
the following result, which was already obtained in \cite{MR3007953,MR1131891}
in a different framework.

\begin{customthm}{\ref{thm:harmonic-stable}}
A harmonic map is always a weakly stable triharmonic map.
\end{customthm}

Later, we will mostly be concerned with a special class of triharmonic maps
which are triharmonic hypersurfaces in space forms. This particular class
of triharmonic maps is explicit enough in order to compute or at least estimate their normal
stability. By normal stability we refer to the second variation of the trienergy
evaluated on vector fields which are proportional to the normal of the hypersurface.

The most general result we establish in this manuscript is given by the following:
\begin{customthm}{\ref{thm:stab-cmc}}
Any triharmonic hypersurface of constant mean curvature in Euclidean space is weakly normally stable
and any triharmonic hypersurface of constant mean curvature in hyperbolic space is normally stable.
\end{customthm}

Following the terminology used for minimal and biharmonic hypersurfaces
we define the normal index of a proper triharmonic hypersurface to 
be the maximal dimension of any linear subspace 
on which the second variation is negative, that is

\begin{align*}
\operatorname{Ind}^{\rm{nor}}(M\to N)
:=\max\{\dim L, L\subset C^\infty_0(M)\mid \operatorname{Hess}E_3(\phi)(f\nu,f\nu)<0
,~~\forall f\in L\},
\end{align*}

where \(\nu\) represents the unit normal of the hypersurface.

Concerning the normal stability of the small proper triharmonic
hypersphere \(\s^m(1/\sqrt{3})\hookrightarrow\s^{m+1}\)
we obtain the following:

\begin{customthm}{\ref{thm:index-small-hypersphere}}
Let \(\phi\colon \s^m(1/\sqrt{3})\hookrightarrow\s^{m+1}\)
be the small proper triharmonic hypersphere.
The index characterizing its normal stability  is 
\begin{align*}
\operatorname{Ind}^{nor}\big(\s^m(1/\sqrt{3})\hookrightarrow\s^{m+1}\big)=1.
\end{align*}
\end{customthm}

Furthermore, in Subsection \ref{sec:clifford} we make some comments on how to estimate the normal index of the proper triharmonic generalized Clifford torus 
$\phi\colon\s^{p}(R_1)\times\s^{q}(R_2)\to\s^{p+q+1}$,
where \(R_1^2+R_2^2=1\) and \(R_1^2\) is determined as the root
of a certain third order polynomial.

\begin{Bem}
If one studies the normal index of the small 
proper biharmonic hypersphere \(\phi\colon\s^m(1/\sqrt{2})\to\s^{m+1}\) one finds, 
as in the triharmonic case, that its normal index is equal to one, see \cite{MR2135286,MR4216418}.
Hence, one may conjecture that the normal index of the small proper \(r\)-harmonic hypersphere 
\(\phi\colon\s^m(1/\sqrt{r})\to\s^{m+1}\), see \cite{MR4106647,MR3711937} for more details, also has normal index equal to one.
\end{Bem}

Throughout this article we will use the following sign conventions:
For the Riemannian curvature
tensor field we use 
\[R(X,Y)Z=[\nabla_X,\nabla_Y]Z-\nabla_{[X,Y]}Z,\]
where \(X,Y,Z\) are vector fields.

For the connection Laplacian on \(\phi^\ast TN\) we use the geometer's
sign convention, that is we set
\(\bar\Delta:=-\tr_g(\bar\nabla\bar\nabla-\bar\nabla_\nabla)\). 
In particular, this implies that the Laplacian acting on functions
has a positive spectrum.

We will make use of the summation convention
and tacitly sum over repeated indices.

This article is organized as follows:
In Section 2 we calculate, mostly for the sake of completeness, the first and second variation of the trienergy
for maps between Riemannian manifolds and then consider the case
that the target manifold is a Riemannian space form of constant curvature.
The third section then first derives a number of general statements on
the normal stability of triharmonic hypersurfaces.
Finally, we calculate the normal index of the small proper triharmonic hypersphere and make some
comments on the normal index of the proper triharmonic Clifford torus.

\section{Variational formulas}
In this section we recall the first and second variation formula
for the trienergy \eqref{3-energy}. To this end we consider a variation of 
\(\phi\), that is \(\phi_t\colon (-\epsilon,\epsilon)\times M\to N\) 
for some small \(\epsilon>0\)
with variational vector field
\begin{align}
\label{dfn:variation-phi}
\frac{\partial\phi_t}{\partial t}\big|_{t=0}=V,
\end{align}
where \(V\in\Gamma(\phi^\ast TN)\).
For simplicity, we will always assume that the variational vector field \(V\)
is compactly supported such that we can employ integration by parts.
We are choosing a local orthonormal frame field \(\{e_i\},i=1,\ldots,\dim M\) tangent to \(M\)
around an arbitrary point \(p\in M\) such that at the point \(p\) we have
\begin{align*}
\nabla_{e_i}e_j=0,i,i=1,\ldots,\dim M,\qquad \nabla_{\partial_t}e_i=0,i=1,\ldots,\dim M.
\end{align*}

In the following we will often make use of the following well-known
formula
\begin{align}
\label{commutator-t-tension}
\frac{\bar\nabla}{\partial t}\tau(\phi_t)=&-\bar\Delta d\phi_t(\partial_t)
+R^N(d\phi_t(\partial_t),d\phi_t(e_i))d\phi_t(e_i).
\end{align}

First of all, we derive the explicit form of the tritension field \(\tau_3(\phi)\).

\begin{Prop}[First Variation]
\label{prop:first-variation}
The critical points of the trienergy \eqref{3-energy} are given by
\begin{align}
0=\tau_3(\phi):=-\bar\Delta^2\tau(\phi)+ R^N(\bar\nabla_{e_j}\tau(\phi),\tau(\phi))d\phi(e_j)
+R^N(\bar\Delta\tau(\phi),d\phi(e_j))d\phi(e_j),
\end{align}
where \(\{e_j\},j=1,\ldots,\dim M\) is a local orthonormal frame field tangent to \(M\).
\end{Prop}

\begin{proof}
We consider a variation of the map \(\phi\) as defined in \eqref{dfn:variation-phi}
and using \eqref{commutator-t-tension} we find
\begin{align*}
\frac{d}{dt}E_3(\phi_t)=&
\int_M\langle\frac{\bar\nabla}{\partial t}\bar\nabla_{e_j}\tau(\phi_t),\bar\nabla_{e_j}\tau(\phi_t)\rangle\dv \\
=&\int_M\big(\langle R^N(d\phi_t(\partial_t),d\phi_t(e_j))\tau(\phi_t),\bar\nabla_{e_j}\tau(\phi_t)\rangle
+\langle \bar\nabla_{e_j}\frac{\bar\nabla}{\partial t}\tau(\phi_t),\bar\nabla_{e_j}\tau(\phi_t) 
\big)\dv
\\
=&\int_M
\big(\langle R^N(\bar\nabla_{e_j}\tau(\phi_t),\tau(\phi_t))d\phi_t(e_j),d\phi_t(\partial_t)\rangle
+\langle \frac{\bar\nabla}{\partial t}\tau(\phi_t),\bar\Delta\tau(\phi_t)\rangle 
\big)\dv \\
=&\int_M
\big(\langle R^N(\bar\nabla_{e_j}\tau(\phi_t),\tau(\phi_t))d\phi_t(e_j),d\phi_t(\partial_t)\rangle
+\langle -\bar\Delta d\phi_t(\partial_t) ,\bar\Delta\tau(\phi_t) \rangle \\
&+\langle R^N(d\phi_t(\partial_t),d\phi_t(e_j))d\phi_t(e_j),\bar\Delta\tau(\phi_t)\rangle
\big)\dv \\
=&\int_M\langle d\phi_t(\partial_t), -\bar\Delta^2\tau(\phi_t)+ R^N(\bar\nabla_{e_j}\tau(\phi_t),\tau(\phi_t))d\phi_t(e_j) \\
&\hspace{1cm}+R^N(\bar\Delta\tau(\phi_t),d\phi_t(e_j))d\phi_t(e_j)\rangle\dv.
\end{align*}
Evaluating at \(t=0\) completes the proof.
\end{proof}

In order to compute the second variation of \eqref{3-energy}
we establish the following commutator formulas.

\begin{Lem}
Consider a variation of the map \(\phi\colon M\to N\) as defined in
\eqref{dfn:variation-phi}. Then, the following formulas hold
\begin{align}
\label{commutator-second-variation}
\frac{\bar\nabla}{\partial t}\bar\nabla_{X}\tau(\phi_t)=&
-\bar\nabla_{X}\bar\Delta d\phi_t(\partial_t)
+\bar\nabla_{X}\big(R^N(d\phi_t(\partial_t),d\phi_t(e_k))d\phi_t(e_k)\big)
+R^N(d\phi_t(\partial_t),d\phi_t(X))\tau(\phi_t), \\
\nonumber\frac{\bar\nabla}{\partial t}\bar\Delta\tau(\phi_t)=&
-R^N(d\phi_t(\partial_t),d\phi_t(e_j))\bar\nabla_{e_j}\tau(\phi_t)
-\bar\nabla_{e_j}\big(R^N(d\phi_t(\partial_t),d\phi_t(e_j))\tau(\phi_t)\big) \\
\nonumber&+\bar\Delta\big(R^N(d\phi_t(\partial_t),d\phi_t(e_j))d\phi_t(e_j)\big)
-\bar\Delta^2d\phi_t(\partial_t), \\
\nonumber\frac{\bar\nabla}{\partial t}\bar\Delta^2\tau(\phi_t)=&
-R^N(d\phi_t(\partial_t),d\phi_t(e_j))\bar\nabla_{e_j}\bar\Delta\tau(\phi_t)
-\bar\nabla_{e_j}\big(R^N(d\phi_t(\partial_t),d\phi_t(e_j))\bar\Delta\tau(\phi_t)\big) \\
&+\nonumber\bar\Delta\frac{\tilde\nabla}{\partial t}\bar\Delta\tau(\phi_t)
\end{align}
for all \(X\in TM\).
\end{Lem}

\begin{proof}
The first equation follows from \eqref{commutator-t-tension} and using
\begin{align*}
\frac{\bar\nabla}{\partial t}\bar\nabla_{X}\tau(\phi_t)
=R^N(d\phi_t(\partial_t),d\phi_t(X))\tau(\phi_t)
+\bar\nabla_{X}\frac{\bar\nabla}{\partial t}\tau(\phi_t),
\end{align*}
where \(X\in TM\).
The second and third equation follow along the same lines.
\end{proof}

\begin{Prop}[Second Variation - rough version]
Let \(\phi\colon M\to N\) be a smooth triharmonic map
and consider a variation of \(\phi\) as defined in \eqref{dfn:variation-phi}.
Then, the second variation of the trienergy is given by
\begin{align}
\label{second-variation-long}
\frac{d^2}{dt^2}\big|_{t=0}E_3(\phi_t)
=\int_M\bigg(&|\bar\nabla\bar\Delta V|^2
+\langle R^N(V,d\phi(e_j))\bar\nabla_{e_j}\bar\Delta\tau(\phi),V\rangle \\
\nonumber&-\langle R^N(V,d\phi(e_j))\bar\Delta\tau(\phi),\bar\nabla_{e_j} V\rangle \\
\nonumber&
+\langle R^N(V,d\phi(e_j))\bar\nabla_{e_j}\tau(\phi),\bar\Delta V\rangle \\
\nonumber&+\langle\bar\nabla_{e_j}\big(R^N(V,d\phi(e_j))\tau(\phi)\big) ,\bar\Delta V\rangle \\
\nonumber&-\langle R^N(V,d\phi(e_j))d\phi(e_j) ,\bar\Delta^2 V\rangle \\
\nonumber&+\langle(\nabla_{V}R^N)(\bar\nabla_{e_j}\tau(\phi),\tau(\phi))d\phi(e_j),V\rangle \\
\nonumber&-\langle R^N(\bar\nabla_{e_j}\bar\Delta V,\tau(\phi))d\phi(e_j),V\rangle \\
\nonumber&+\langle R^N(\bar\nabla_{e_j}\big(R^N(V,d\phi(e_k))d\phi(e_k)\big),\tau(\phi))d\phi(e_j),V\rangle \\
\nonumber&+\langle R^N(R^N(V,d\phi(e_j))\tau(\phi),\tau(\phi))d\phi(e_j),V\rangle \\
\nonumber&-\langle R^N(\bar\nabla_{e_j}\tau(\phi),\bar\Delta V)d\phi(e_j),V\rangle \\
\nonumber&+\langle R^N(\bar\nabla_{e_j}\tau(\phi),R^N(V,d\phi(e_i))d\phi(e_i))d\phi(e_j),V\rangle \\
\nonumber&+\langle R^N(\bar\nabla_{e_j}\tau(\phi),\tau(\phi))\bar\nabla_{e_j}V,V\rangle \\
\nonumber&+\langle(\nabla_{V}R^N)(\bar\Delta\tau(\phi),d\phi(e_j))d\phi(e_j),V\rangle \\
\nonumber&-\langle R^N(R^N(V,d\phi(e_j))\bar\nabla_{e_j}\tau(\phi),d\phi(e_k))d\phi(e_k),V\rangle \\
\nonumber&-\langle R^N(\bar\nabla_{e_j}\big(R^N(V,d\phi(e_j))\tau(\phi)\big),d\phi(e_k))d\phi(e_k),V\rangle \\  
\nonumber&+\langle R^N(\bar\Delta\big(R^N(V,d\phi(e_j))d\phi(e_j)\big),d\phi(e_k))d\phi(e_k),V\rangle \\  
\nonumber&-\langle R^N(\bar\Delta^2V,d\phi(e_j))d\phi(e_j),V\rangle \\  
\nonumber&+\langle R^N(\bar\Delta\tau(\phi),\bar\nabla_{e_j}V)d\phi(e_j),V\rangle \\
\nonumber&+\langle R^N(\bar\Delta\tau(\phi),d\phi(e_j))\bar\nabla_{e_j}V,V\rangle\bigg)\dv.
\end{align}

\end{Prop}

\begin{proof}
We consider a variation of the map \(\phi\) as defined in \eqref{dfn:variation-phi}
and use the first variational formula derived in Proposition \ref{prop:first-variation}
as a starting point. Now, we compute
\begin{align*}
\frac{d^2}{dt^2}E_3(\phi_t)=&
-\int_M\langle\frac{\bar\nabla}{\partial t}\bar\Delta^2\tau(\phi_t),d\phi_t(\partial_t)\rangle\dv \\
&+\int_M\langle\frac{\bar\nabla}{\partial t}\big(R^N(\bar\nabla_{e_j}\tau(\phi_t),\tau(\phi_t))d\phi_t(e_j)\big),d\phi_t(\partial_t)\rangle\dv \\
&+\int_M\langle \frac{\bar\nabla}{\partial t}\big(R^N(\bar\Delta\tau(\phi_t),d\phi_t(e_j))d\phi_t(e_j)\big),
d\phi_t(\partial_t)\rangle \dv
+\int_M\langle\underbrace{\tau_3(\phi_t)}_{=0},\frac{\bar\nabla}{\partial t}d\phi_t(\partial_t)\rangle\dv \\
:=&\int_M\big(A_1+A_2+A_3\big)\dv.
\end{align*}

We will now calculate all three terms individually.
The \(A_1\) contribution can be manipulated as follows
\begin{align*}
A_1=&\langle\frac{\bar\nabla}{\partial t}\bar\Delta^2\tau(\phi_t),d\phi_t(\partial_t)\rangle \\
=&-\langle R^N(d\phi_t(\partial_t),d\phi_t(e_j))\bar\nabla_{e_j}\bar\Delta\tau(\phi_t),d\phi_t(\partial_t)\rangle \\
&-\langle\bar\nabla_{e_j}\big(R^N(d\phi_t(\partial_t),d\phi_t(e_j))\bar\Delta\tau(\phi_t)\big),
d\phi_t(\partial_t)\rangle
+\langle\bar\Delta\frac{\tilde\nabla}{\partial t}\bar\Delta\tau(\phi_t),d\phi_t(\partial_t)\rangle,
\end{align*}
where we used the third equation of \eqref{commutator-second-variation}.

In addition, we get
\begin{align*}
\langle\frac{\tilde\nabla}{\partial t}\bar\Delta\tau(\phi_t),\bar\Delta d\phi_t(\partial_t)\rangle
=&
-\langle R^N(d\phi_t(\partial_t),d\phi_t(e_j))\bar\nabla_{e_j}\tau(\phi_t),\bar\Delta d\phi_t(\partial_t)\rangle \\
&-\langle\bar\nabla_{e_j}\big(R^N(d\phi_t(\partial_t),d\phi_t(e_j))\tau(\phi_t)\big),
\bar\Delta d\phi_t(\partial_t)\rangle \\
&+\langle\bar\Delta\big(R^N(d\phi_t(\partial_t),d\phi_t(e_j))d\phi_t(e_j)\big),\bar\Delta d\phi_t(\partial_t)\rangle \\
&-\langle\bar\Delta^2d\phi_t(\partial_t) ,\bar\Delta d\phi_t(\partial_t)\rangle
\end{align*}
using the second equation of \eqref{commutator-second-variation}.

The second term, \(A_2\), can be rewritten as follows
\begin{align*}
A_2=&\langle\frac{\bar\nabla}{\partial t}\big(R^N(\bar\nabla_{e_j}\tau(\phi_t),
\tau(\phi_t))d\phi_t(e_j)\big),d\phi_t(\partial_t)\rangle \\
=&\langle(\nabla_{d\phi_t(\partial_t)}R^N)(\bar\nabla_{e_j}\tau(\phi_t),\tau(\phi_t))d\phi_t(e_j),d\phi_t(\partial_t)\rangle \\
&+\langle R^N(\frac{\bar\nabla}{\partial t}\bar\nabla_{e_j}\tau(\phi_t),\tau(\phi_t))d\phi_t(e_j),d\phi_t(\partial_t)\rangle \\
&+\langle R^N(\bar\nabla_{e_j}\tau(\phi_t),\frac{\bar\nabla}{\partial t}\tau(\phi_t))d\phi_t(e_j),d\phi_t(\partial_t)\rangle \\
&+\langle R^N(\bar\nabla_{e_j}\tau(\phi_t),
\tau(\phi_t))\bar\nabla_{e_j}d\phi_t(\partial_t),d\phi_t(\partial_t)\rangle\\
=&\langle(\nabla_{d\phi_t(\partial_t)}R^N)(\bar\nabla_{e_j}\tau(\phi_t),\tau(\phi_t))d\phi_t(e_j),d\phi_t(\partial_t)\rangle \\
&-\langle R^N(\bar\nabla_{e_j}\bar\Delta d\phi_t(\partial_t),\tau(\phi_t))d\phi_t(e_j),d\phi_t(\partial_t)\rangle \\
&+\langle R^N(\bar\nabla_{e_j}\big(R^N(d\phi_t(\partial_t),d\phi_t(e_k))d\phi_t(e_k)\big),\tau(\phi_t))d\phi_t(e_j),d\phi_t(\partial_t)\rangle \\
&+\langle R^N(R^N(d\phi_t(\partial_t),d\phi_t(e_j))\tau(\phi_t),\tau(\phi_t))d\phi_t(e_j),d\phi_t(\partial_t)\rangle \\
&-\langle R^N(\bar\nabla_{e_j}\tau(\phi_t),\bar\Delta d\phi_t(\partial_t))d\phi_t(e_j),d\phi_t(\partial_t)\rangle \\
&+\langle R^N(\bar\nabla_{e_j}\tau(\phi_t),R^N(d\phi_t(\partial_t),d\phi_t(e_i))d\phi_t(e_i))d\phi_t(e_j),d\phi_t(\partial_t)\rangle\\
&+\langle R^N(\bar\nabla_{e_j}\tau(\phi_t),
\tau(\phi_t))\bar\nabla_{e_j}d\phi_t(\partial_t),d\phi_t(\partial_t)\rangle,
\end{align*}
where we used the first equation of \eqref{commutator-second-variation} after the second
equals sign and \eqref{commutator-t-tension} in the final step.

Concerning the \(A_3\) term we find
\begin{align*}
A_3=&\langle\frac{\bar\nabla}{\partial t}\big(R^N(\bar\Delta\tau(\phi_t),d\phi_t(e_j))d\phi_t(e_j)\big),d\phi_t(\partial_t)\rangle\\
=&\langle(\nabla_{d\phi_t(\partial_t)}R^N)(\bar\Delta\tau(\phi_t),d\phi_t(e_j))d\phi_t(e_j),d\phi_t(\partial_t)\rangle \\
&+\langle R^N(\frac{\bar\nabla}{\partial t}\bar\Delta\tau(\phi_t),d\phi_t(e_j))d\phi_t(e_j),d\phi_t(\partial_t)\rangle  \\
&+\langle R^N(\bar\Delta\tau(\phi_t),\bar\nabla_{e_j}d\phi_t(\partial_t))d\phi_t(e_j),d\phi_t(\partial_t)\rangle \\
&+\langle R^N(\bar\Delta\tau(\phi_t),d\phi_t(e_j))\bar\nabla_{e_j}d\phi_t(\partial_t),d\phi_t(\partial_t)\rangle \\
=&\langle(\nabla_{d\phi_t(\partial_t)}R^N)(\bar\Delta\tau(\phi_t),d\phi_t(e_j))d\phi_t(e_j),d\phi_t(\partial_t)\rangle \\
&-\langle R^N(R^N(d\phi_t(\partial_t),d\phi_t(e_j))\bar\nabla_{e_j}\tau(\phi_t),d\phi_t(e_k))d\phi_t(e_k),d\phi_t(\partial_t)\rangle \\
&-\langle R^N(\bar\nabla_{e_j}\big(R^N(d\phi_t(\partial_t),d\phi_t(e_j))\tau(\phi_t)\big),d\phi_t(e_k))d\phi_t(e_k),d\phi_t(\partial_t)\rangle \\  
&+\langle R^N(\bar\Delta\big(R^N(d\phi_t(\partial_t),d\phi_t(e_j))d\phi_t(e_j)\big),d\phi_t(e_k))d\phi_t(e_k),d\phi_t(\partial_t)\rangle \\  
&-\langle R^N(\bar\Delta^2d\phi_t(\partial_t),d\phi_t(e_j))d\phi_t(e_j),d\phi_t(\partial_t)\rangle \\  
&+\langle R^N(\bar\Delta\tau(\phi_t),\bar\nabla_{e_j}d\phi_t(\partial_t))d\phi_t(e_j),d\phi_t(\partial_t)\rangle \\
&+\langle R^N(\bar\Delta\tau(\phi_t),d\phi_t(e_j))\bar\nabla_{e_j}d\phi_t(\partial_t),d\phi_t(\partial_t)\rangle,
\end{align*}
where we employed the second equation of \eqref{commutator-second-variation} in the second step.
Adding up the different contributions and evaluating at \(t=0\)
then completes the proof.
\end{proof}

Luckily, many of the terms in \eqref{second-variation-long}
can be grouped together such that
we can obtain the following simplification:

\begin{Prop}[Second Variation]
Let \(\phi\colon M\to N\) be a smooth triharmonic map
and consider a variation of \(\phi\) as defined in \eqref{dfn:variation-phi}
with \(V\) compactly supported.
Then, the second variation of the trienergy \eqref{3-energy} is given by
\begin{align}
\label{second-variation-simplified}
\frac{d^2}{dt^2}\big|_{t=0}E_3(\phi_t)
=\int_M\bigg(&\big|\bar\nabla\bar\Delta V-
\bar\nabla\big(R^N(V,d\phi(e_k))d\phi(e_k)\big)\big|^2 
+|R^N(V,d\phi(e_j))\tau(\phi)|^2\\
\nonumber&-2\langle R^N(V,d\phi(e_j))\bar\Delta\tau(\phi),\bar\nabla_{e_j} V\rangle\\
\nonumber&+2\langle R^N(V,d\phi(e_j))\bar\nabla_{e_j}\tau(\phi),\bar\Delta V\rangle,\\
\nonumber&-2\langle R^N(V,d\phi(e_k))d\phi(e_k),\bar\nabla_{e_j}\big(R^N(V,d\phi(e_j))\tau(\phi)\big)\rangle\\
\nonumber&+2\langle R^N(d\phi(e_j),V)\bar\nabla_{e_j}\tau(\phi),R^N(V,d\phi(e_k))d\phi(e_k)\rangle \\
\nonumber&-2\langle R^N(\bar\nabla_{e_j}\bar\Delta V,\tau(\phi))d\phi(e_j),V\rangle\\
\nonumber&+\langle(\nabla_{V}R^N)(\bar\nabla_{e_j}\tau(\phi),\tau(\phi))d\phi(e_j),V\rangle \\
\nonumber&+\langle(\nabla_{V}R^N)(\bar\Delta\tau(\phi),d\phi(e_j))d\phi(e_j),V\rangle \\
\nonumber&+\langle R^N(V,d\phi(e_j))\bar\nabla_{e_j}\bar\Delta\tau(\phi),V\rangle \\
\nonumber&+\langle R^N(\bar\Delta\tau(\phi),d\phi(e_j))\bar\nabla_{e_j}V,V\rangle \\
\nonumber&+\langle R^N(\bar\nabla_{e_j}\tau(\phi),\tau(\phi))\bar\nabla_{e_j}V,V\rangle 
\bigg)\dv.
\end{align}
\end{Prop}

\begin{proof}
First of all, we note that
\begin{align*}
\int_M\langle\bar\nabla_{e_j}\big(R^N(V,d\phi(e_j))\tau(\phi)\big) ,\bar\Delta V\rangle\dv
=-\int_M\langle R^N(V,d\phi(e_j))\tau(\phi),\bar\nabla_{e_j}\bar\Delta V\rangle\dv
\end{align*}
such that
\begin{align*}
\int_M&\big(\langle\bar\nabla_{e_j}\big(R^N(V,d\phi(e_j))\tau(\phi)\big),\bar\Delta V\rangle
-\langle R^N(\bar\nabla_{e_j}\bar\Delta V,\tau(\phi))d\phi(e_j),V\rangle
\big)\dv\\
=&-2\int_M\langle R^N(\bar\nabla_{e_j}\bar\Delta V,\tau(\phi))d\phi(e_j),V\rangle\dv.
\end{align*}

Using the symmetries of the Riemann curvature tensor we get
\begin{align*}
-\langle R^N&(V,d\phi(e_j))\bar\Delta\tau(\phi),\bar\nabla_{e_j} V\rangle
+\langle R^N(\bar\Delta\tau(\phi),\bar\nabla_{e_j}V)d\phi(e_j),V\rangle \\
&=-2\langle R^N(V,d\phi(e_j))\bar\Delta\tau(\phi),\bar\nabla_{e_j} V\rangle,\\
\langle R^N&(V,d\phi(e_j))\bar\nabla_{e_j}\tau(\phi),\bar\Delta V\rangle
-\langle R^N(\bar\nabla_{e_j}\tau(\phi),\bar\Delta V)d\phi(e_j),V\rangle \\
&=2\langle R^N(V,d\phi(e_j))\bar\nabla_{e_j}\tau(\phi),\bar\Delta V\rangle.
\end{align*}

Moreover, we find
\begin{align*}
&\langle R^N(\bar\nabla_{e_j}\tau(\phi),R^N(V,d\phi(e_i))d\phi(e_i))d\phi(e_j),V\rangle
-\langle R^N(R^N(V,d\phi(e_j))\bar\nabla_{e_j}\tau(\phi),d\phi(e_k))d\phi(e_k),V\rangle \\
&=2\langle R^N(d\phi(e_j),V)\bar\nabla_{e_j}\tau(\phi),R^N(V,d\phi(e_k))d\phi(e_k)\rangle.
\end{align*}

Once more using the symmetries of the Riemann curvature tensor we obtain
\begin{align*}
\langle R^N(R^N(V,d\phi(e_j))\tau(\phi),\tau(\phi))d\phi(e_j),V\rangle
=|R^N(V,d\phi(e_j))\tau(\phi)|^2.
\end{align*}

By a similar reasoning we get
\begin{align*}
-\langle R^N(\bar\nabla_{e_j}\big(R^N(V,d\phi(e_j))\tau(\phi)\big),d\phi(e_k))d\phi(e_k),V\rangle \\
=-\langle R^N(V,d\phi(e_k))d\phi(e_k),\bar\nabla_{e_j}\big(R^N(V,d\phi(e_j))\tau(\phi)\big)\rangle
\end{align*}
and also
\begin{align*}
\langle R^N(\bar\nabla_{e_j}\big(R^N(V,d\phi(e_k))d\phi(e_k)\big),\tau(\phi))d\phi(e_j),V\rangle \\
=\langle R^N(V,d\phi(e_j))\tau(\phi),\bar\nabla_{e_j}\big(R^N(V,d\phi(e_k))d\phi(e_k)\big)\rangle.
\end{align*}
Hence, using integration by parts we can infer
\begin{align*}
\int_M\big(&-\langle R^N(\bar\nabla_{e_j}\big(R^N(V,d\phi(e_j))\tau(\phi)\big),d\phi(e_k))d\phi(e_k),V\rangle \\
&+\langle R^N(\bar\nabla_{e_j}\big(R^N(V,d\phi(e_k))d\phi(e_k)\big),\tau(\phi))d\phi(e_j),V\rangle\big)\dv \\
&=-2\int_M\langle R^N(V,d\phi(e_k))d\phi(e_k),\bar\nabla_{e_j}\big(R^N(V,d\phi(e_j))\tau(\phi)\big)\rangle\dv.
\end{align*}
Using both the symmetries of the Riemann curvature tensor and integration by parts we can
deduce
\begin{align*}
\int_M&\langle R^N(\bar\Delta\big(R^N(V,d\phi(e_j))d\phi(e_j)\big),d\phi(e_k))d\phi(e_k),V\rangle\dv \\
&=\int_M\langle R^N(V,d\phi(e_j))d\phi(e_j),\bar\Delta\big(R^N(V,d\phi(e_k))d\phi(e_k)\big)\rangle\dv \\
&=\int_M|\bar\nabla\big(R^N(V,d\phi(e_j))d\phi(e_j)\big)|^2\dv.
\end{align*}

Finally, we note that
\begin{align*}
-\langle R^N&(V,d\phi(e_j))d\phi(e_j) ,\bar\Delta^2 V\rangle
-\langle R^N(\bar\Delta^2V,d\phi(e_j))d\phi(e_j),V\rangle \\
&=-2\langle R^N(V,d\phi(e_j))d\phi(e_j) ,\bar\Delta^2 V\rangle
\end{align*}
which allows to complete the square in the first term of \eqref{second-variation-simplified}
finishing the proof.
\end{proof}

A direct application of the previous Proposition is the following 
result which was already given in \cite[Proposition 3.4]{MR3007953}, see also \cite{MR1131891},
using the formula for the second variation of polyharmonic maps.
\begin{Theorem}
\label{thm:harmonic-stable}
A harmonic map is always a weakly stable triharmonic map.
\end{Theorem}
\begin{proof}
By assumption we have \(\tau(\phi)=0\) such that \eqref{second-variation-simplified}
simplifies to
\begin{align*}
\label{}
\frac{d^2}{dt^2}\big|_{t=0}E_3(\phi_t)
=\int_M\big(&\big|\bar\nabla\bar\Delta V-
\bar\nabla\big(R^N(V,d\phi(e_k))d\phi(e_k)\big)\big|^2 
\big)\dv\geq 0
\end{align*}
which already yields the result.
\end{proof}

\subsection{The second variation formula for triharmonic maps to space forms}
In order to further simplify the formula for the second variation \eqref{second-variation-simplified} we will now consider the case that the target manifold has a
particular simple geometric structure.
More precisely, we will study the stability of triharmonic maps
in the case that the ambient space is a space form of constant curvature \(K\).
Then, the Riemann curvature tensor acquires the simple form
\begin{align}
\label{curvature-space-form}
R^N(X,Y)Z=K(\langle Y,Z\rangle X-\langle X,Z\rangle Y),
\end{align}
where \(X,Y,Z\) are vector fields on \(N\) and \(K\)
represents the constant curvature of the space form \(N\).

In the following we will use the notation 
\begin{align*}
\operatorname{Hess}E_3(\phi)(V,V):=\frac{d^2}{dt^2}\big|_{t=0}E_3(\phi_t)
\end{align*}
for a variation of \(\phi\) that satisfies \(\frac{\bar\nabla\phi_t}{\partial t}\big|_{t=0}=V\)
still assuming that the variational vector field \(V\) is compactly supported.

\begin{Prop}
Let \(\phi\colon M\to N\) be a smooth triharmonic map and assume that \(N\)
is a space form of constant curvature \(K\).
Then, the second variation of the trienergy \eqref{3-energy} simplifies to
\begin{align}
\label{second-variation-space-form}
\operatorname{Hess}E_3(\phi)(V,V)
=\int_M\bigg(&
\big|\bar\nabla\bar\Delta V-
\bar\nabla\big(R^N(V,d\phi(e_k))d\phi(e_k)\big)\big|^2\\
\nonumber&
+K\big[
3\langle V,\bar\Delta\tau(\phi)\rangle\langle d\phi,\bar\nabla V\rangle 
-2\langle d\phi(e_j),\bar\Delta\tau(\phi)\rangle\langle V,\bar\nabla_{e_j}V\rangle 
\\
\nonumber&+2\langle d\phi,\bar\nabla\tau(\phi)\rangle\langle V,\bar\Delta V\rangle
-2\langle V,\bar\nabla_{e_j}\tau(\phi)\rangle\langle d\phi(e_j),\bar\Delta V\rangle \\
\nonumber&
-2\langle\tau(\phi),d\phi(e_j)\rangle\langle\bar\nabla_{e_j}\bar\Delta V,V\rangle
+2\langle\bar\nabla\bar\Delta V,d\phi\rangle\langle\tau(\phi),V\rangle \\
\nonumber&
+\langle d\phi,\bar\nabla\bar\Delta\tau(\phi)\rangle|V|^2
-\langle V,\bar\nabla_{e_j}\bar\Delta\tau(\phi)\rangle\langle V,d\phi(e_j)\rangle \\
\nonumber&
+\langle\tau(\phi),\bar\nabla_{e_j}V\rangle\langle\bar\nabla_{e_j}\tau(\phi),V\rangle
-\langle \bar\nabla\tau(\phi),\bar\nabla V\rangle\langle V,\tau(\phi)\rangle \\
\nonumber&-\langle\bar\Delta\tau(\phi),\bar\nabla_{e_j}V\rangle\langle d\phi(e_j),V\rangle 
\big]\\
\nonumber&+K^2\big[
3|d\phi|^2|\langle V,\tau(\phi)\rangle|^2
+|\langle d\phi,\tau(\phi)\rangle|^2|V|^2 
-2|d\phi|^2|\tau(\phi)|^2|V|^2 \\
\nonumber&+4|d\phi|^2\langle V,\bar\nabla_{e_i}\tau(\phi)\rangle\langle V,d\phi(e_i)\rangle
\\
\nonumber&+4|\langle d\phi,V\rangle|^2\langle d\phi,\bar\nabla\tau(\phi)\rangle
-4|d\phi|^2|V|^2\langle d\phi,\bar\nabla\tau(\phi)\rangle \\
\nonumber&-2\langle d\phi(e_i),\tau(\phi)\rangle\langle V,\tau(\phi)\rangle
\langle V,d\phi(e_i)\rangle
-2|d\phi|^2\langle d\phi(e_i),\tau(\phi)\rangle\langle V,\bar\nabla_{e_i} V\rangle \\
\nonumber&+2|d\phi|^2\langle V,d\phi(e_i)\rangle\langle\bar\nabla_{e_i} V,\tau(\phi)\rangle
+2|\tau(\phi)|^2|\langle d\phi,V\rangle|^2 \\
\nonumber&-4\langle V,\bar\nabla_{e_j}\tau(\phi)\rangle
\langle d\phi(e_j),d\phi(e_k)\rangle\langle d\phi(e_k),V\rangle \\
\nonumber&+2\langle d\phi(e_k),V\rangle\langle d\phi(e_j),\tau(\phi)\rangle 
 \langle d\phi(e_k),\bar\nabla_{e_j}V\rangle \\
 \nonumber&-2\langle d\phi(e_k),V\rangle\langle\bar\nabla_{e_j}V,\tau(\phi)\rangle
 \langle d\phi(e_k),d\phi(e_j)\rangle \\
 \nonumber&-2\langle d\phi(e_k),V\rangle\langle d\phi(e_k),\tau(\phi)\rangle
 \langle V,\tau(\phi)\rangle
\big]
 \bigg)\dv.
\end{align}
\end{Prop}

\begin{proof}
We make use of \eqref{second-variation-simplified} and the fact that the target manifold \(N\) is a space form of constant curvature \(K\). Due to this assumption the terms involving the derivative of the curvature tensor in \eqref{second-variation-simplified} drop out.

Now, using \eqref{curvature-space-form} we find
\begin{align*}
\langle R^N(V,d\phi(e_j))\bar\Delta\tau(\phi),\bar\nabla_{e_j} V\rangle 
&=K\big(\langle d\phi(e_j),\bar\Delta\tau(\phi)\rangle\langle V,\bar\nabla_{e_j}V\rangle 
-\langle V,\bar\Delta\tau(\phi)\rangle\langle d\phi,\bar\nabla V\rangle\big),\\
\langle R^N(V,d\phi(e_j))\bar\nabla_{e_j}\tau(\phi),\bar\Delta V\rangle
&=K\big(\langle d\phi,\bar\nabla\tau(\phi)\rangle\langle V,\bar\Delta V\rangle
-\langle V,\bar\nabla_{e_j}\tau(\phi)\rangle\langle d\phi(e_j),\bar\Delta V\rangle.
\end{align*}

In addition, we get
\begin{align*}
\langle R^N(\bar\nabla_{e_j}\bar\Delta V,\tau(\phi))d\phi(e_j),V\rangle
=&K\big(\langle\tau(\phi),d\phi(e_j)\rangle\langle\bar\nabla_{e_j}\bar\Delta V,V\rangle
-\langle\bar\nabla\bar\Delta V,d\phi\rangle\langle\tau(\phi),V\rangle\big),\\
\langle R^N(V,d\phi(e_j))\bar\nabla_{e_j}\bar\Delta\tau(\phi),V\rangle 
=&K\big(\langle d\phi,\bar\nabla\bar\Delta\tau(\phi)\rangle|V|^2
-\langle V,\bar\nabla_{e_j}\bar\Delta\tau(\phi)\rangle\langle V,d\phi(e_j)\rangle\big), \\
\langle R^N(\bar\Delta\tau(\phi),d\phi(e_j))\bar\nabla_{e_j}V,V\rangle
=&K\big(\langle d\phi,\bar\nabla V\rangle\langle\bar\Delta\tau(\phi),V\rangle
-\langle\bar\Delta\tau(\phi),\bar\nabla_{e_j}V\rangle\langle d\phi(e_j),V\rangle\big), \\
\langle R^N(\bar\nabla_{e_j}\tau(\phi),\tau(\phi))\bar\nabla_{e_j}V,V\rangle 
=&K\big(\langle\tau(\phi),\bar\nabla_{e_j}V\rangle\langle\bar\nabla_{e_j}\tau(\phi),V\rangle
-\langle \bar\nabla\tau(\phi),\bar\nabla V\rangle\langle V,\tau(\phi)\rangle\big).
\end{align*}

Moreover, we have
\begin{align*}
\langle R^N(d\phi(e_j),V)&\bar\nabla_{e_j}\tau(\phi),R^N(V,d\phi(e_k))d\phi(e_k)\rangle \\
=&K^2\big(|d\phi|^2\langle V,\bar\nabla_{e_j}\tau(\phi)\rangle\langle d\phi(e_j),V\rangle
-\langle V,\bar\nabla_{e_j}\tau(\phi)\rangle\langle d\phi(e_j),d\phi(e_k)\rangle\langle d\phi(e_k),V\rangle \\
&-|d\phi|^2|V|^2\langle d\phi,\bar\nabla\tau(\phi)\rangle 
+|\langle d\phi,V\rangle|^2\langle d\phi,\bar\nabla\tau(\phi)\rangle\big). 
\end{align*}

A similar analysis yields
\begin{align*}
|R^N(V,d\phi(e_j))\tau(\phi)|^2
=K^2\big(&|\langle d\phi,\tau(\phi)\rangle|^2|V|^2+|d\phi|^2|\langle V,\tau(\phi)\rangle|^2 \\
&-2\langle d\phi(e_i),\tau(\phi)\rangle\langle V,\tau(\phi)\rangle
\langle V,d\phi(e_i)\rangle\big).
\end{align*}

In order to rewrite the term \(\langle R^N(V,d\phi(e_k))d\phi(e_k),\bar\nabla_{e_j}\big(R^N(V,d\phi(e_j))\tau(\phi)\big)\rangle
\)
we first calculate
\begin{align*}
\bar\nabla_{e_j}\big(R^N(V,d\phi(e_j))\tau(\phi)\big)
=&K\big(|\tau(\phi)|^2V+\langle d\phi,\bar\nabla\tau(\phi)\rangle V
+\langle d\phi(e_j),\tau(\phi)\rangle\bar\nabla_{e_j}V\\
&-\langle \bar\nabla_{e_j}V,\tau(\phi)\rangle d\phi(e_j)
-\langle V,\bar\nabla_{e_j}\tau(\phi)\rangle d\phi(e_j)
-\langle V,\tau(\phi)\rangle\tau(\phi)\big).
\end{align*}
Now, we can infer that
\begin{align*}
\langle& R^N(V,d\phi(e_k))d\phi(e_k),\bar\nabla_{e_j}\big(R^N(V,d\phi(e_j))\tau(\phi)\big)\rangle \\
=&K^2\big(
|d\phi|^2|\tau(\phi)|^2|V|^2+|d\phi|^2|V|^2\langle d\phi,\bar\nabla\tau(\phi)\rangle
+|d\phi|^2\langle d\phi(e_j),\tau(\phi)\rangle\langle V,\bar\nabla_{e_j}V\rangle \\
&-|d\phi|^2\langle V,d\phi(e_j)\rangle\langle\bar\nabla_{e_j}V,\tau(\phi)\rangle
-|d\phi|^2\langle V,\bar\nabla_{e_j}\tau(\phi)\rangle\langle V,d\phi(e_j)\rangle
-|d\phi|^2|\langle V,\tau(\phi)\rangle|^2 \\
&-|\tau(\phi)|^2|\langle d\phi,V\rangle|^2
-|\langle d\phi,V\rangle|^2\langle d\phi,\bar\nabla\tau(\phi)\rangle
-\langle d\phi(e_k),V\rangle\langle d\phi(e_j),\tau(\phi)\rangle
\langle d\phi(e_k),\bar\nabla_{e_j}V\rangle \\
&+\langle d\phi(e_k),V\rangle\langle\bar\nabla_{e_j}V,\tau(\phi)\rangle\langle d\phi(e_k),d\phi(e_j)\rangle
+\langle d\phi(e_k),V\rangle\langle d\phi(e_k),d\phi(e_j)\rangle\langle V,\bar\nabla_{e_j}\tau(\phi)\rangle \\
&+\langle d\phi(e_k),V\rangle\langle d\phi(e_k),\tau(\phi)\rangle\langle V,\tau(\phi)\rangle\big).
\end{align*}
Adding up the different contributions now completes the proof.
\end{proof}

\begin{Bem}
We could of course also rewrite the first term in \eqref{second-variation-space-form}
using the fact that the target manifold is a space form of constant curvature \(K\).
However, as this term already has a fixed sign it seems most convenient to
not further manipulate it at this stage of the computation.
\end{Bem}

In the following we will mostly be concerned with the normal stability of triharmonic
hypersurfaces in space forms by choosing a variational vector field \(V\)
that takes values in the normal bundle of the hypersurface.
Under this assumption we have \(\langle V,d\phi\rangle=0\)
such that we get the following simplification of the previous Lemma.

\begin{Cor}
\label{cor:normal-space-form}
Let \(\phi\colon M\to N\) be a smooth triharmonic map and assume that \(N\)
is a space form of constant curvature \(K\).
If we assume that \(\langle V,d\phi\rangle=0\)
the second variation of the trienergy \eqref{3-energy} simplifies to
\begin{align}
\label{second-variation-space-form-normal}
\operatorname{Hess}E_3(\phi)(V,V)
=\int_M\bigg(&
\big|\bar\nabla\bar\Delta V-
K\bar\nabla(|d\phi|^2V)\big|^2\\
\nonumber&
+K\big[
3\langle V,\bar\Delta\tau(\phi)\rangle\langle d\phi,\bar\nabla V\rangle 
-2\langle d\phi(e_j),\bar\Delta\tau(\phi)\rangle\langle V,\bar\nabla_{e_j}V\rangle \\
\nonumber&+2\langle d\phi,\bar\nabla\tau(\phi)\rangle\langle V,\bar\Delta V\rangle
-2\langle V,\bar\nabla_{e_j}\tau(\phi)\rangle\langle d\phi(e_j),\bar\Delta V\rangle \\
\nonumber&
-2\langle\tau(\phi),d\phi(e_j)\rangle\langle\bar\nabla_{e_j}\bar\Delta V,V\rangle
+2\langle\bar\nabla\bar\Delta V,d\phi\rangle\langle\tau(\phi),V\rangle \\
\nonumber&
+\langle\tau(\phi),\bar\nabla_{e_j}V\rangle\langle\bar\nabla_{e_j}\tau(\phi),V\rangle
-\langle \bar\nabla\tau(\phi),\bar\nabla V\rangle\langle V,\tau(\phi)\rangle \\
\nonumber&
+\langle d\phi,\bar\nabla\bar\Delta\tau(\phi)\rangle|V|^2 
\big]\\
\nonumber&+K^2\big[
3|d\phi|^2|\langle V,\tau(\phi)\rangle|^2
+|\langle d\phi,\tau(\phi)\rangle|^2|V|^2 
-2|d\phi|^2|\tau(\phi)|^2|V|^2 \\
\nonumber&
-4|d\phi|^2|V|^2\langle d\phi,\bar\nabla\tau(\phi)\rangle 
-2|d\phi|^2\langle d\phi(e_j),\tau(\phi)\rangle\langle V,\bar\nabla_{e_j} V\rangle\big]
 \bigg)\dv.
    \end{align}
\end{Cor}
\begin{proof}
This follows directly from \eqref{second-variation-space-form} 
using that \(\langle d\phi,V\rangle=0\).
\end{proof}

\section{The normal stability of triharmonic hypersurfaces in space forms}
In this section we apply the formula for the second variation of the trienergy 
obtained in the previous section in order to study the normal stability
of a particular class of triharmonic maps.

\subsection{The normal stability of CMC triharmonic hypersurfaces}
First, we will derive a number of general statements on the normal stability of proper triharmonic hypersurfaces.
We will only consider the CMC (constant mean curvature) case as this already comes with many computational difficulties and as all known, at least up to now, triharmonic hypersurfaces
belong to this category.

In order to perform the computations outlined above, let us recall
a number of geometric facts on hypersurfaces \(M^m\) in a Riemannian manifold 
\(N^{m+1}\).
The connections on \(M^m\) and \(N^{m+1}\) are related by the following formula
\begin{align}
\label{dfn:sff}
\nabla^M_XY=\nabla^N_XY+\sff(X,Y),
\end{align}
where \(\sff\) represents the second fundamental form of the hypersurface.

Let \(\nu\) be the global unit normal of the hypersurface \(M^m\), then 
its shape operator \(A\) is given by
\begin{align}
\label{dfn:shape}
\nabla_{X}\nu=-A(X),
\end{align}
where \(X\) is a vector field on \(M\).

The shape operator \(A\) and the second fundamental form \(\sff\)
are related by
\begin{align}
\label{relation-shape-sff}
\sff(X,Y)=\langle A(X),Y\rangle\nu.
\end{align}

If \(\phi\colon M^{m}\to N^{m+1}\) is an isometric immersion
the tension field acquires the form \(\tau(\phi)=mH\nu\),
where \(H\) represents the mean curvature function of the hypersurface.
In addition, we can compute the mean curvature function by \(H=\frac{1}{m}\tr A\).

In the following lemma we provide without proof three standard facts which we shall use in this section, 
see for example \cite[Lemma 2.1]{MR4462636}
or \cite[Lemma 4.1]{MR4552081}.

\begin{Lem}
\label{lemma-sff}
Let $\phi\colon M^m\to N^{m+1}$ be a hypersurface in a Riemannian manifold \(N^{m+1}\).
In addition, let $A$ be the shape operator and $H=1/m\tr A$ the mean curvature function. Then, we have that 
\begin{enumerate}
\item $(\nabla A) (\cdot,\cdot)$ is symmetric,
\item $\langle (\nabla A) (\cdot,\cdot), \cdot \rangle$ is totally symmetric,
\item $\tr (\nabla A) (\cdot,\cdot)= m\operatorname{grad} H$. 
\end{enumerate}
\end{Lem}

Now, we establish the following Lemma,
the first statement was already established in \cite[Lemma 2.2]{MR4462636}.

\begin{Lem}
Let $\phi\colon M^m\to N^{m+1}$ be a hypersurface. 
Suppose that \(k\in C^\infty(M)\) and let \(\nu\) be the unit normal
of the hypersurface. Then, the following identities hold
\begin{align}
\label{eq:laplace-hypersurface}
\bar\Delta(k\nu)=&(\Delta k+k|A|^2)\nu+2A(\operatorname{grad} k)+mk\operatorname{grad} H,\\
\label{eq:nabla-laplace-hypersurface}
\bar\nabla_{X}\bar\Delta(k\nu)=&
(\nabla_{X}\Delta k+(\nabla_{X} k)|A|^2+k\nabla_{X}|A|^2)\nu 
-(\Delta k+k|A|^2)A(X) \\
\nonumber&+2(\nabla_{X}A)(\operatorname{grad} k)
+2A(\nabla_{X}\operatorname{grad} k)
+2\sff(X,A(\operatorname{grad} k))\\
&\nonumber+m(\nabla_{X}k)\operatorname{grad} H
+mk\nabla_{X}\operatorname{grad} H
+mk\sff(X,\operatorname{grad} H),
\end{align}
where \(X\) represents a vector field on \(M\).
\end{Lem}

\begin{proof}
Here, we closely follow the calculations carried out in \cite[Lemma 2.2]{MR4462636}
using a slightly different notation.
We choose a local orthonormal frame \(\{e_i\},i=1,\ldots,m\)
that satisfies \(\nabla_{e_j}e_i=0\) at a fixed point \(p\in M\)
for all \(i,j=1,\ldots,m\).
Now, we compute
\begin{align*}
\bar\nabla_{e_i}(k\nu)=(\nabla_{e_i}k)\nu+k\nabla_{e_i}\nu=(\nabla_{e_i}k)\nu-kA(e_i)
\end{align*}
and differentiating once more we get
\begin{align*}
\bar\nabla_{e_i}\bar\nabla_{e_i}(k\nu)=&
(\nabla_{e_i}\nabla_{e_i}k)\nu+\nabla_{e_i}k\nabla_{e_i}\nu
-\nabla_{e_i}kA(e_i)-k\big(\nabla A(e_i,e_i)+\sff(e_i,A(e_i))\big)\\
=&(\nabla_{e_i}\nabla_{e_i}k)\nu-2\nabla_{e_i}kA(e_i)-k(\nabla A(e_i,e_i))
-k|A(e_i)|^2\nu.
\end{align*}
Now, the first claim follows from summing over \(i\) and using Lemma \ref{lemma-sff}.

Concerning the second equation \eqref{eq:nabla-laplace-hypersurface}, 
we use the first equation, differentiate with respect to \(X\)
and apply both \eqref{dfn:sff} and \eqref{dfn:shape}.
\end{proof}

In order to study the normal stability of triharmonic hypersurfaces in space forms
of constant curvature \(K\),
we will now employ the general formula \eqref{second-variation-space-form}
and choose the section \(V=f\nu\) with \(f\in C_0^\infty(M)\).
Note that this choice of \(V\) trivializes the normal bundle of the hypersurface.

In the following we will make use of the  
identity
\begin{align}
\label{eq:identity-square-variation}
\int_M&|\bar\nabla\bar\Delta(f\nu)-mK\bar\nabla(f\nu)|^2\dv \\
\nonumber&=\int_M|\bar\nabla\bar\Delta(f\nu)|^2-2mK|\bar\Delta(f\nu)|^2
+m^2K^2\big(|\nabla f|^2+f^2|A|^2\big)\dv
\end{align}
which follows from integration by parts and
the definition of the shape operator.

\begin{Prop}
Let \(\phi\colon M^m\to N^{m+1}\) be a smooth CMC triharmonic hypersurface
with unit normal \(\nu\)
and assume that \(N\) is a space form of constant curvature \(K\).
Then, the second variation of the trienergy evaluated at \(V=f\nu\)
is given by
\begin{align}
\label{second-variation-hypersurface-cmc}
\operatorname{Hess}E_3(\phi)(f\nu,f\nu)
=\int_M\big(&|\bar\nabla\bar\Delta(f\nu)|^2-2mK|\bar\Delta(f\nu)|^2
+2mHK\langle\bar\nabla\bar\Delta V,d\phi\rangle f
\\
\nonumber&+m^2K^2\big(|\nabla f|^2+|A|^2f^2\big) \\
\nonumber&+5m^3K^2H^2f^2
-7m^2K|A|^2 H^2f^2
-2m^2KH^2 f\Delta f
\big)\dv,
\end{align}
where \(f\in C_0^\infty(M)\).
\end{Prop}

\begin{proof}
In order to prove the claim we have to evaluate all contributions 
in \eqref{second-variation-space-form-normal} in the geometric setup of a
hypersurface with constant mean curvature, that is
\(\tau(\phi)=mH\nu\) with \(H=const\) and insert the variational vector field 
\(V=f\nu\).
We will explain this in detail for the third and fourth term of \eqref{second-variation-space-form-normal} and just state the final result for the other contributions.
As before, we choose a local orthonormal frame \(\{e_i\},i=1,\ldots,m\)
that satisfies \(\nabla_{e_j}e_i=0\) at a fixed point \(p\in M\)
for all \(i,j=1,\ldots,m\).

Using \eqref{eq:laplace-hypersurface} with \(k=1\) we find
\begin{align*}
\bar\Delta\tau(\phi)=&mH\Delta\nu=mH|A|^2\nu,
\end{align*}
where we also made use of the assumption \(H=const\). Hence, we get that
\(\langle d\phi(X),\bar\Delta\tau(\phi)\rangle=0\) 
for all vector fields \(X\)
and thus
\begin{align*}
\langle d\phi(e_j),\bar\Delta\tau(\phi)\rangle\langle V,\bar\nabla_{e_j}V\rangle=0.
\end{align*}

In addition, we find
\begin{align*}
\underbrace{\langle V,\bar\Delta\tau(\phi)\rangle}_{=mH|A|^2f}
\underbrace{\langle d\phi,\bar\nabla V\rangle}
_{=-\langle\tau(\phi),V\rangle=-mHf}=-m^2|A|^2 H^2f^2.
\end{align*}

Regarding the other terms, for \(V=f\nu\), similar computations lead to
\begin{align*}
\langle d\phi,\bar\nabla\tau(\phi)\rangle\langle V,\bar\Delta V\rangle=& 
-m^2H^2f(\Delta f+f|A|^2),\\
\langle V,\bar\nabla_{e_j}\tau(\phi)\rangle\langle d\phi(e_j),\bar\Delta V\rangle=&0, 
\\
\langle d\phi,\bar\nabla\bar\Delta\tau(\phi)\rangle|V|^2=&-m^2|A|^2H^2f^2, \\
\langle\tau(\phi),\bar\nabla_{e_j}V\rangle\langle\bar\nabla_{e_j}\tau(\phi),V\rangle=& 
0,\\ 
\langle\bar\nabla\tau(\phi),\bar\nabla V\rangle\langle V,\tau(\phi)\rangle=&
m^2H^2|A|^2f^2,
\\
|d\phi|^2|V|^2\langle d\phi,\bar\nabla\tau(\phi)\rangle=&-m^3H^2f^2,\\
\langle\tau(\phi),d\phi(e_j)\rangle\langle\bar\nabla_{e_j}\bar\Delta V,V\rangle=&0,\\
\langle\bar\nabla\bar\Delta V,d\phi\rangle\langle\tau(\phi),V\rangle=&
mfH\langle\bar\nabla\bar\Delta V,d\phi\rangle.
\end{align*}
Replacing all terms in \eqref{second-variation-space-form-normal}
and using the identity \eqref{eq:identity-square-variation}
then completes the proof.
\end{proof}

An immediate consequence of the above formula are the following facts:
\begin{Prop}
Any minimal hypersurface
in a space form of constant curvature \(K\)
is a weakly stable triharmonic hypersurface with respect to normal variations.
\end{Prop}
\begin{proof}
Using \eqref{eq:identity-square-variation}
in \eqref{second-variation-hypersurface-cmc} 
together with the assumption \(H=0\)
we find
\begin{align}
\operatorname{Hess}E_3(\phi)(f\nu,f\nu)
=\int_M\big(|\bar\nabla\bar\Delta(f\nu)-mK\bar\nabla(f\nu)|^2\big)\dv\geq 0,
\end{align}
which already completes the proof.
\end{proof}

\begin{Thm}
\label{thm:stab-cmc}
Any triharmonic hypersurface of constant mean curvature in Euclidean space
is weakly normally stable 
and any triharmonic hypersurface of constant mean curvature in hyperbolic space is normally stable.
\end{Thm}
\begin{proof}
In order to obtain the claim we first derive the estimate
\begin{align*}
2mHK\langle\bar\nabla\bar\Delta V,d\phi\rangle f\geq &
-2m|H||K||\bar\nabla\bar\Delta V||d\phi||f| \\
\geq &-|\bar\nabla\bar\Delta V|^2-m^2H^2K^2|d\phi|^2f^2\\
= &-|\bar\nabla\bar\Delta V|^2-m^3H^2K^2f^2,
\end{align*}
where we employed Young's inequality.

Inserting into \eqref{second-variation-hypersurface-cmc} we find
\begin{align*}
\operatorname{Hess}E_3(\phi)(f\nu,f\nu)
&=\int_M\big(-2mK|\bar\Delta(f\nu)|^2
+m^2K^2\big(|\nabla f|^2+|A|^2f^2\big) \\
\nonumber&\hspace{1cm}+4m^3K^2H^2f^2
-7m^2K|A|^2 H^2f^2
-2m^2KH^2 |\nabla f|^2
\big)\dv \\
&\geq 0,
\end{align*}
where made use of the assumption \(K\leq 0\) in the last step 
completing the proof.
\end{proof}

We conclude that the only interesting case to study in detail the normal stability
of a triharmonic hypersurface in space forms of constant curvature
is in the case of a spherical target.
In order to systematically approach this question we establish the following
\begin{Lem}
Let \(\phi\colon M^m\to N^{m+1}\)
be a hypersurface with constant mean curvature and \(|A|^2=const\).
Let \(\nu\) be the unit normal of the hypersurface and \(f\in C^\infty_0(M)\).
Then the following formulas hold
\begin{align}
\label{eq:identity-a}
|\bar\Delta (f\nu)|^2
% &(\Delta f+f|A|^2)^2+\big(2A(\operatorname{grad} f)\big)^2\\
=&|\Delta f|^2+f^2|A|^4+2\Delta f f|A|^2 
+4|A(\operatorname{grad} f)|^2, \\
\label{eq:identity-b}
|\bar\nabla\bar\Delta(f\nu)|^2 
=&|\nabla\Delta f|^2+|\nabla f|^2|A|^4+2\nabla\Delta f\nabla f|A|^2 \\
\nonumber&+|\Delta f|^2|A|^2+f^2|A|^6+2f\Delta f|A|^4 \\
\nonumber&+4|(\nabla A)(\operatorname{grad} f)|^2
+4|A(\nabla\operatorname{grad} f)|^2
\nonumber+4|\sff(e_i,A(\operatorname{grad} f))|^2\\
\nonumber&+8\langle(\nabla A)(\operatorname{grad} f),A(\nabla\operatorname{grad} f)\rangle
+8\langle(\nabla_{e_i} A)(\operatorname{grad} f),\sff(e_i,A(\operatorname{grad} f))\rangle \\
\nonumber&+4\nabla\Delta f\langle(\nabla A)(\operatorname{grad} f),\nu\rangle
+4|A|^2\nabla f\langle(\nabla A)(\operatorname{grad} f),\nu\rangle \\
\nonumber&+4\langle\sff(\operatorname{grad}\Delta f,A(\operatorname{grad} f)),\nu\rangle
+4|A|^2\langle\sff(\operatorname{grad} f,A(\operatorname{grad} f)),\nu\rangle \\
\nonumber&-4\Delta fA((\nabla A)(\operatorname{grad}f))
-4f|A|^2A((\nabla A)(\operatorname{grad}f)) \\
\nonumber&-4\Delta f\langle A(e_i),A(\nabla_{e_i}\operatorname{grad} f)\rangle
-4f|A|^2\langle A(e_i),A(\nabla_{e_i}\operatorname{grad} f)\rangle.
\end{align}
\end{Lem}

\begin{proof}
Using the assumption that the hypersurface has constant mean curvature
we find
\begin{align*}
\bar\Delta (f\nu)=
(\Delta f+f|A|^2)\nu+2A(\operatorname{grad} f),
\end{align*}
the first claim now follows from computing the square.

In order to prove the second claim we employ \eqref{eq:nabla-laplace-hypersurface}
with \(|A|^2=const\) and \(H=const\).
Then, for all vector fields \(X\), we find
\begin{align*}
\bar\nabla_{X}\bar\Delta(f\nu)=&
(\nabla_{X}\Delta f+(\nabla_{X} f)|A|^2)\nu 
-(\Delta f+f|A|^2)A(X) \\
\nonumber&
+2(\nabla_{X}A)(\operatorname{grad} f)
+2A(\nabla_{X}\operatorname{grad} f)
+2\sff(X,A(\operatorname{grad} f)),
\end{align*}
where we used the assumption that the hypersurface has constant mean curvature.
The claim now follows from a direct, but lengthy computation.
\end{proof}

After having established the necessary geometric formulas for triharmonic 
hypersurfaces in Riemannian manifolds 
we also recall a number of results on the spectrum of the Laplace operator on closed manifolds.
On an arbitrary closed Riemannian manifold \((M,g)\)
its eigenvalues are discrete, have finite multiplicity and satisfy
\begin{align*}
0=\lambda_0\leq\lambda_1\leq\lambda_2\leq\ldots\to\infty.
\end{align*}

In the case of the \(p\)-dimensional sphere of radius \(R\), that is \(\s^p(R)\),
we know that the eigenvalues of the Laplacian are given by
\begin{align}
\label{ev-laplacian-sphere}
\lambda_j=\frac{j(j+p-1)}{R^2},\qquad j=0,1,2,\ldots
\end{align}
with multiplicities \(m_{\lambda_0}=1,m_{\lambda_1}=p+1\)
and
\begin{align*}
m_{\lambda_j}={{p+j}\choose{j}}-{{p+j-2}\choose{j-2}},\qquad j=2,\ldots
\end{align*}

Now, let \((M_i,g_i),i=1,2\) be two closed Riemannian manifolds.
Then, we know that if \(\lambda\) is an eigenvalue of the Laplacian on \((M_1,g_1)\)
with multiplicity \(m_\lambda\) and \(\mu\) an eigenvalue
of the Laplacian on \((M_2,g_2)\) with multiplicity \(m_\mu\),
then \(\nu=\lambda+\mu\) is an eigenvalue of the Laplacian on \(M_1\times M_2\).
The multiplicity of the eigenvalue \(\nu\) is the sum of the products 
\(m_\lambda m_\mu\) of all \(\lambda\) and \(\mu\) which satisfy
\(\nu=\lambda+\mu\).

For the general properties of the spectrum of the Laplacian we refer to the book \cite{MR0282313}, for more details concerning the Laplace operator in the context of
the stability of CMC hypersurfaces one may consult
\cite[Section 3]{MR2336585}.

\subsection{The normal stability of the small triharmonic hypersphere}
As a next step we will investigate the normal stability of proper triharmonic hypersurfaces
in the sphere. There are two well-known examples of such proper triharmonic hypersurfaces
\cite{MR4106647,MR3711937},
the small hypersphere \(\phi\colon\s^m(1/\sqrt{3})\hookrightarrow\s^{m+1}\)
and a certain generalized Clifford torus 
\(\phi\colon\s^{p}(R_1) \times \s^{q}(R_2)\hookrightarrow\s^{p+q+1}\) with \(R_1^2+R_2^2=1\),
see Subsection \ref{sec:clifford} for the precise definition.

The following result from \cite[p. 356]{MR4462636} characterizes proper triharmonic hypersurfaces
which have constant norm of the shape operator.
\begin{Prop}
\label{prop:mor-tri-hyper}
A hypersurface \(\phi\colon M^{m}\to N^{m+1}\), where \(N^{m+1}\) is
a space form of constant curvature \(K\), is triharmonic if
\begin{align*}
\Delta|A|^2+|A|^4-mK|A|^2-m^2KH^2=0,\qquad A(\operatorname{grad}|A|^2)=0.
\end{align*}
Here, \(|A|^2\) represents the norm of the shape operator.
\end{Prop}

Hence, in the case that \(|A|^2\) is constant, the algebraic condition for the hypersurface
being triharmonic is
\begin{align*}
|A|^4-mK|A|^2-m^2KH^2=0.
\end{align*}
It becomes clear that in order to find proper triharmonic hypersurfaces
the only interesting case is \(K=1\) as
for \(K=-1\) we would have that \(|A|^2=0\) corresponding to a minimal hypersurface.

Moreover, let us also recall the following fact:
Consider the inclusion \(\phi\colon\s^m(a)\hookrightarrow\s^{m+1}\),
then the corresponding shape operator is given by
\begin{align}
\label{shape-operator-inclusion}
A=-\frac{\sqrt{1-a^2}}{a}\operatorname{Id}.
\end{align}

Hence, concerning the small hypersphere \(\phi\colon\s^m(1/\sqrt{3})\hookrightarrow\s^{m+1}\)
we thus have \(|A|^2=2m\) and \(H^2=2\).

With these geometric data we are ready to give the following
\begin{Prop}
Let \(\phi\) be the proper triharmonic hypersphere \(\phi\colon\s^m(1/\sqrt{3})\hookrightarrow
\s^{m+1}\).
The quadratic form describing its normal stability is given by
\begin{align}
\label{normal-stability-small-first}
Q(V,V)=
\int_{\s^m(\frac{1}{\sqrt{3}})}&\big(
|\nabla\Delta f|^2
+(4m+16)|\Delta f|^2 
+8|\nabla^2 f|^2
+(-3m^2+16+20m)|\nabla f|^2
\\
\nonumber&
-24m^3f^2\big)\dv,
\end{align}
where we chose the variational vector field
\(V=f\nu\) with \(f\in C_0^\infty(M)\) and \(\nu\)
represents the unit normal of the hypersurface.
\end{Prop}

\begin{proof}
We employ \eqref{second-variation-hypersurface-cmc}, set \(K=1\)
and use the identities \eqref{eq:identity-a}, \eqref{eq:identity-b}.
By \eqref{shape-operator-inclusion} we have that \(\nabla A=0\)
and \(|A|^2=2m\). Now, a lengthy, but direct calculation yields
\begin{align*}
|\bar\nabla\bar\Delta(f\nu)|^2-2m|\bar\Delta (f\nu)|^2
=&|\nabla\Delta f|^2+4m^2|\nabla f|^2+4m\nabla\Delta f\nabla f \\
&+4|A(\nabla\operatorname{grad})f|^2
+4|\sff(e_i,A(\operatorname{grad} f))|^2\\
&+4\langle\sff(\operatorname{grad}\Delta f,A(\operatorname{grad} f)),\nu\rangle
+8m\langle\sff(\operatorname{grad} f,A(\operatorname{grad} f)),\nu\rangle \\
&-4\Delta f\langle A(e_i),A(\nabla_{e_i}\operatorname{grad} f)\rangle 
-8mf\langle A(e_i),A(\nabla_{e_i}\operatorname{grad} f)\rangle \\
&-8m|A(\operatorname{grad} f)|^2.
\end{align*}

As a next step we replace all terms involving the second fundamental form
by the shape operator using \eqref{relation-shape-sff}, that is
\begin{align*}
|\sff(e_i,A(\operatorname{grad} f))|^2
&=|\langle A(e_i),A(\operatorname{grad} f)\rangle|^2
,\\
\langle\sff(\operatorname{grad}\Delta f,A(\operatorname{grad} f)),\nu\rangle
&=\langle A(\operatorname{grad}\Delta f),A(\operatorname{grad} f)\rangle
,\\
\langle\sff(\operatorname{grad} f,A(\operatorname{grad} f)),\nu\rangle
&=\langle A(\operatorname{grad} f),A(\operatorname{grad} f)\rangle.
\end{align*}

Now, by setting \(a=\frac{1}{\sqrt{3}}\) in \eqref{shape-operator-inclusion} 
we have for all vector fields \(X\) the identity
\begin{align*}
A(X)=-\sqrt{2}X
\end{align*}
such that we get
\begin{align*}
|\bar\nabla\bar\Delta(f\nu)|^2-2m|\bar\Delta (f\nu)|^2
=&|\nabla\Delta f|^2+4m^2|\nabla f|^2+4m\nabla\Delta f\nabla f+8|\nabla^2 f|^2
+16|\nabla f|^2\\
&+8\nabla\Delta f\nabla f
+8|\Delta f|^2+16mf\Delta f.
\end{align*}
Furthermore, we employ integration by parts
\begin{align*}
\int_{\s^m(\frac{1}{\sqrt{3}})}\langle\bar\nabla\bar\Delta (f\nu),d\phi\rangle f\dv
=&-\int_{\s^m(\frac{1}{\sqrt{3}})}\langle\bar\Delta (f\nu),\tau(\phi)\rangle f\dv
-\int_{\s^m(\frac{1}{\sqrt{3}})}\langle\bar\Delta (f\nu),d\phi(e_j)\rangle\nabla_{e_j}f\dv \\
=&-\int_{\s^m(\frac{1}{\sqrt{3}})}(\Delta f+|A|^2f)mHf\dv \\
&-2\int_{\s^m(\frac{1}{\sqrt{3}})}\langle A(\operatorname{grad}f),d\phi(e_j)\rangle
\nabla_{e_j}f\dv,
\end{align*}
where we used \eqref{eq:laplace-hypersurface} in the second step.

Using that \(A(X)=-\sqrt{2}X\) for all vector fields \(X\) we can infer
\begin{align*}
\langle A(\operatorname{grad}f),d\phi(e_j)\rangle
\nabla_{e_j}f=-\sqrt{2}|\nabla f|^2.
\end{align*}
Hence, we may conclude that
\begin{align*}
2mH\int_{\s^m(\frac{1}{\sqrt{3}})}\langle\bar\nabla\bar\Delta (f\nu),d\phi\rangle f\dv
=4(-m^2+m)\int_{\s^m(\frac{1}{\sqrt{3}})}|\nabla f|^2\dv
-8m^3\int_{\s^m(\frac{1}{\sqrt{3}})}f^2\dv.
\end{align*}
Finally, replacing \(H^2=2\) and \(|A|^2=2m\) once more, we obtain
for the remaining terms of \eqref{second-variation-hypersurface-cmc}
that
\begin{align*}
\int_{\s^m(\frac{1}{\sqrt{3}})}&\big(m^2\big(|\nabla f|^2+f^2|A|^2\big)
+5m^3f^2H^2 
-7m^2|A|^2 H^2f^2
-2m^2H^2 f\Delta f
\big)\dv \\
&=-\int_{\s^m(\frac{1}{\sqrt{3}})}(3m^2|\nabla f|^2+16m^3f^2)\dv
\end{align*}
and combining the previous equations completes the proof.
\end{proof}

After having fixed the geometric setup for the stability analysis of triharmonic
hypersurfaces in the sphere we now make the following definition.

\begin{Dfn}
Let \(\phi\colon M^m\to N^{m+1}\) be a proper triharmonic hypersurface.
Then, its normal index is defined to be the maximal dimension 
of any linear subspace \(L\subset C^\infty_0(M)\)
on which \(\operatorname{Hess}E_3(\phi)(f\nu,f\nu)\) is negative, that is
\begin{align}
\label{definition-normal-index}
\operatorname{Ind}^{\rm{nor}}(M)
:=\max\{\dim L, L\subset C^\infty_0(M)\mid \operatorname{Hess}E_3(\phi)(f\nu,f\nu)<0
,~~\forall f\in L\}.
\end{align}
\end{Dfn}

\begin{Bem}
We can already read of from \eqref{normal-stability-small-first}
that the index of the small proper triharmonic hypersphere \(\phi\colon\s^m(1/\sqrt{3})\to\s^{m+1}\)
is at least one: If we choose \(f=const\) corresponding to the eigenfunction of the Laplacian
associated with the eigenvalue zero, the quadratic form \eqref{normal-stability-small-first} is negative.
\end{Bem}

In the following, we will make use of the classic Bochner formula
\begin{align*}
\int_M|\nabla^2f|^2\dv=-\int_M\operatorname{Ric}^M(\nabla f,\nabla f)\dv
+\int_M|\Delta f|^2\dv,
\end{align*}
which holds on every closed Riemannian manifolds.

In the case that \(M=\s^m(a)\)
we have \(\operatorname{Ric}^M=\frac{m-1}{a^2}g\) such that
the above formula simplifies to
\begin{align}
\label{bochner-function-sphere}
\int_{\s^m(\frac{1}{\sqrt{3}})}|\nabla^2f|^2\dv
=-3(m-1)\int_{\s^m(\frac{1}{\sqrt{3}})}|\nabla f|^2\dv
+\int_{\s^m(\frac{1}{\sqrt{3}})}|\Delta f|^2\dv.
\end{align}

With these preparations we are now ready to give the following
\begin{Theorem}
\label{thm:index-small-hypersphere}
The normal index of the small proper triharmonic hypersphere 
\(\phi\colon\s^m(1/\sqrt{3})\hookrightarrow\s^{m+1}\)
is equal to one, that is
\begin{align*}
\operatorname{Ind}^{\rm{nor}}(\s^m(1/\sqrt{3})\hookrightarrow\s^{m+1})=1.
\end{align*}
\end{Theorem}

\begin{proof}
In order to obtain the claim
we assume that \(f\) is an eigenfunction of the Laplacian, that is \(\Delta f=\lambda f\).
Now, together with \eqref{bochner-function-sphere} we then get from \eqref{normal-stability-small-first} that
\begin{align*}
Q(V,V)=
\int_{\s^m(\frac{1}{\sqrt{3}})}
\big(\lambda^3
+4(m+6)\lambda^2
+(-3m^2-4m+40)\lambda
-24 m^3
\big)f^2\dv.
\end{align*}
Recall that the first non-zero eigenvalue of the Laplace operator on \(\s^m(\frac{1}{\sqrt{3}})\) is given by \(\lambda_1=3m\), see \eqref{ev-laplacian-sphere}, in which case we find
\begin{align*}
Q(V,V)=
\int_{\s^m(\frac{1}{\sqrt{3}})}
\big(30m^3+204m^2+120m
\big)f^2\dv>0.
\end{align*}
Now, it is easy to see that the above quadratic form will be positive for
all higher eigenvalues of the Laplace operator.
As the dimension of the vector space of constant functions, for which the quadratic form \(Q(V,V)\)
is negative, is one-dimensional we have finished the proof.
\end{proof}

\subsection{The normal stability of the generalized triharmonic Clifford torus}
\label{sec:clifford}
There is a second explicit example of a CMC proper triharmonic hypersurface
in the sphere which is given by a certain generalized Clifford torus.

Before turning to the proper triharmonic Clifford torus let us recall
the following facts which hold true for every Clifford torus
$\phi\colon\s^{p}(R_1)\times\s^{q}(R_2)\hookrightarrow\s^{p+q+1}$
with \(R_1^2+R_2^2=1\).
For \(X_1\in\Gamma(T\s^p(R_1))\) and \(X_2\in\Gamma(T\s^q(R_2))\)
the corresponding shape operator satisfies
\begin{align}
\label{shape-clifford}
A(X_1)=-\frac{R_2}{R_1}X_1,\qquad A(X_2)=\frac{R_1}{R_2}X_2
\end{align}
such that the following identities hold
\begin{align*}
H=\frac{1}{p+q}\big(-\frac{R_2}{R_1}p+\frac{R_1}{R_2}q\big)\nu,\qquad
|A|^2=\frac{R_2^2}{R_1^2}p+\frac{R_1^2}{R_2^2}q.
\end{align*}

For more details on the stability of both minimal and CMC Clifford tori we refer to
\cite{MR2301375,MR2336585}.

In order to characterize proper triharmonic Clifford tori we recall
the following 
\begin{Theorem}[\cite{MR4106647,MR3711937}]
Let $p,q \geq 1 $ and assume that the radii $R_1,R_2$ satisfy $R_1^2+R_2^2=1$.
Then, the generalized Clifford torus $\phi\colon\s^{p}(R_1)\times\s^{q}(R_2)\hookrightarrow\s^{p+q+1}$ is proper triharmonic if it is not minimal and if 
$x=R_1^2$ is a root of the following polynomial
\begin{align}
\label{triharmonic-clifford-polynomial}
P(x)=3(p+q)x^3-(2q+7p)x^2+5px-p.
\end{align}
\end{Theorem}

It is straightforward to check that \(q=0\) and \(x=\frac{1}{3}\) is a root
of \eqref{triharmonic-clifford-polynomial} which corresponds to
the small hypersphere studied in the previous section.

For \(p,q\geq 1\), using some facts on cubic equations it is straightforward to see
that \eqref{triharmonic-clifford-polynomial} has, in general, at least one real root.

In addition, it can be checked directly that in the case \(p=q\) the only root of
\eqref{triharmonic-clifford-polynomial} is \(x=\frac{1}{2}\) corresponding
to a minimal Clifford torus, see \cite[Remark 2.13]{MR4106647}. 
We already know from Theorem \ref{thm:stab-cmc}
that such hypersurfaces are normally stable such that we do not further need to investigate their 
normal stability.

Luckily, we also have a second way to describe the above hypersurface.
In the case that \(|A|^2\) is constant, we know from Proposition \ref{prop:mor-tri-hyper} 
that the condition for being a proper triharmonic hypersurface in
the sphere is
\begin{align}
\label{shape-clifford-algebra}
|A|^4-m|A|^2-m^2H^2=0
\end{align}
which we will employ in the following.

\begin{Bem}
Since the shape operator of the generalized Clifford torus \eqref{shape-clifford} 
is no longer proportional to the identity it becomes substantially more difficult
to investigate its normal stability
compared to the case of the small hypersphere.
Note that in the case of the biharmonic generalized Clifford torus this problem does not occur
as the corresponding index, describing the normal stability, only contains 
the norm of the shape operator but not the shape operator itself, see \cite[Theorem 3.1]{MR4386842}.
For these reasons it seems that we are only able to derive an upper bound
on the normal index of the triharmonic generalized Clifford torus.  
\end{Bem}

First, we will prove the following estimate which holds
for every proper triharmonic hypersurface in the sphere.
\begin{Prop} 
Let \(\phi\colon M^m\to\s^{m+1}\) be a proper CMC triharmonic hypersurface
with constant norm of the shape operator.
The quadratic form describing its normal stability 
can be estimated as
\begin{align}
\label{stability-general}
\operatorname{Hess}E_3(\phi)(f\nu,f\nu)
\geq\int_{M}\bigg(&
|\nabla\Delta f|^2 
-2m|\Delta f|^2 -4|A|^2|\nabla f||\nabla\Delta f|\\ 
\nonumber&-2\sqrt{m}|A|^2\sqrt{|A|^2-m}|f||\Delta f|
+(m^2-10m|A|^2)|\nabla f|^2 \\
\nonumber&+|A|^2(7m|A|^2-m^2
-2\sqrt{m}|A|^2\sqrt{|A|^2-m}f^2
-9|A|^4)
\bigg)\dv,
\end{align}
where \(f\in C_0^\infty(M)\) and \(\nu\)
represents the normal of the hypersurface.
\end{Prop}

\begin{proof}
From \eqref{shape-clifford} we have that \(\nabla A=0\), replacing all terms involving
the second fundamental form using \eqref{relation-shape-sff},
equations \eqref{eq:identity-a} and \eqref{eq:identity-b} yield 
\begin{align*}
\int_M\bigg(&|\bar\nabla\bar\Delta(f\nu)|^2-2m |\bar\Delta (f\nu)|^2
+2mH\langle\bar\nabla\bar\Delta(f\nu),d\phi\rangle\big)\dv \\
=&\int_M\big(|\nabla\Delta f|^2
+(3|A|^2-2m)|\Delta f|^2
+(3|A|^4-4m|A|^2)|\nabla f|^2 
+|A|^4(|A|^2-2m)f^2 
\\
&-8m|A(\operatorname{grad} f)|^2
+4|A(\nabla\operatorname{grad} f)|^2
+4|\langle A(e_i),A(\operatorname{grad} f)\rangle|^2\\
&+4\langle A(\nabla\Delta f),A(\operatorname{grad} f)\rangle
+4|A|^2\langle A(\nabla f),A(\operatorname{grad} f)\rangle \\
&-4\Delta f\langle A(e_i),A(\nabla_{e_i}\operatorname{grad} f)\rangle
-4f|A|^2\langle A(e_i),A(\nabla_{e_i}\operatorname{grad} f)\rangle \\
&-2mH \langle d\phi(e_i),A(e_i)\rangle f\Delta f
-2mH|A|^2\langle d\phi(e_i),A(e_i)\rangle f^2 \\
&+4mH\langle A(\nabla_{e_i}\operatorname{grad} f),d\phi(e_i)\rangle f
\bigg)\dv,
\end{align*}
where we also used integration by parts.

As a next step we estimate all terms involving the shape operator
such that we get
\begin{align*}
\int_M\big(&|\bar\nabla\bar\Delta(f\nu)|^2-2m |\bar\Delta (f\nu)|^2
+2mH\langle\bar\nabla\bar\Delta(f\nu),d\phi\rangle
\big)\dv \\
\geq &
\int_M\bigg(|\nabla\Delta f|^2
+(3|A|^2-2m)|\Delta f|^2
+(3|A|^4-12m|A|^2)|\nabla f|^2 
+|A|^4(|A|^2-2m)f^2 \\
&+4|A(\nabla\operatorname{grad} f)|^2
+4|\langle A(e_i),A(\operatorname{grad} f)\rangle|^2\\
&-4|\langle A(e_i),A(\operatorname{grad} f)\rangle|\big(|\nabla\Delta f|+|A|^2|\nabla f|\big) \\
&-4|A(\nabla\operatorname{grad} f)|\big(|\Delta f||A|+|f||A|^3+m|H||d\phi||f|\big) \\
&-2m|H||d\phi||A||f||\Delta f|
-2m|H||A|^3|d\phi|f^2
\bigg)\dv.
\end{align*}

Note that we have the following estimates 
\begin{align*}
|A(\nabla\operatorname{grad}& f)|^2
-|A(\nabla\operatorname{grad} f)|\big(|A||\Delta f|+|A|^3f+m|H||d\phi||f|\big) \\
&\geq-\frac{1}{4}\big(|A||\Delta f|+|A|^3|f|+m|H||d\phi||f|\big)^2\\
&=-\frac{3}{4}|A|^2|\Delta f|^2-\frac{3}{4}|A|^6f^2-\frac{3}{4}m^3H^2f^2,
\end{align*}
where we made use of Young's inequality.

Once more, we estimate
\begin{align*}
|\langle A(e_i)&,A(\operatorname{grad} f)\rangle|^2
-|\langle A(e_i),A(\operatorname{grad} f)\rangle|\big(|\nabla\Delta f|+|A|^2|\nabla f|\big) \\
&\geq
-|A|^2|\nabla f||\nabla\Delta f|-\frac{1}{4}|A|^4|\nabla f|^2,
\end{align*}
where we again used Young's inequality.
From \eqref{shape-clifford-algebra} we get \(H=|A|\frac{\sqrt{|A|^2-m}}{m}\) such that
\begin{align*}
-2m&|H||d\phi||A||f||\Delta f|
-2m|H||A|^3|d\phi|f^2 \\
&=-2\sqrt{m}|A|^2\sqrt{|A|^2-m}|f||\Delta f|
-2\sqrt{m}|A|^4\sqrt{|A|^2-m}f^2.
\end{align*}
By combining these estimates we arrive at 
\begin{align*}
\int_M\big(&|\bar\nabla\bar\Delta(f\nu)|^2-2m |\bar\Delta (f\nu)|^2+2mH\langle\bar\nabla\bar\Delta(f\nu),d\phi\rangle\big)\dv \\
\geq&\int_M\big(
|\nabla\Delta f|^2
-2m|\Delta f|^2
+2|A|^2(|A|^2-6m)|\nabla f|^2 
+|A|^2(3m^2-5m|A|^2-2|A|^4)f^2 
\\
&-4|A|^2|\nabla\Delta f||\nabla f|
-2\sqrt{m}|A|^2\sqrt{|A|^2-m}|f||\Delta f|
-2\sqrt{m}|A|^4\sqrt{|A|^2-m}f^2
\big)\dv.
\end{align*}
The remaining terms in \eqref{second-variation-hypersurface-cmc}
can easily be manipulated by using \(H^2=\frac{|A|^4}{m^2}-\frac{|A|^2}{m}\) such that we get
\begin{align*}
\int_{M}&\big(m^2|\nabla f|^2+m^2f^2|A|^2+5m^3f^2H^2-7m^2|A|^2 H^2f^2-2m^2H^2 f\Delta f\big)\dv \\
=&\int_{M}\big(
(m^2-2|A|^4+2m|A|^2)|\nabla f|^2+(-4m^2|A|^2+12m|A|^4-7|A|^6)f^2
\big)\dv.
\end{align*}
Combining all terms now completes the proof.
\end{proof}

Using the estimate \eqref{stability-general} we can, at least in principal,
now give a bound on the normal index of the proper triharmonic Clifford torus
by performing the following steps:

\begin{enumerate}
 \item Fix the dimensions \(p\) and \(q\) and determine \(R_1^2\) via \eqref{triharmonic-clifford-polynomial}.
 \item Using the value of \(R_1^2\), calculate \(|A|^2\) via \eqref{shape-clifford} and 
 the eigenvalues of the Laplace operator as explained after equation \eqref{ev-laplacian-sphere}.
 \item Insert this data into \eqref{stability-general} and determine the first eigenvalue 
 of the Laplace operator for which \eqref{stability-general} becomes positive. 
 The multiplicity of this eigenvalue then determines the normal index of the proper triharmonic Clifford torus.
\end{enumerate}

Although the above algorithm can easily be implemented into a computer algebra system
we do not present any further details as the estimate provided
by \eqref{stability-general} does not seem to be very accurate.

\bibliographystyle{plain}
\bibliography{mybib}
\end{document}